\documentclass{amsart}
\usepackage{amssymb, amsmath, color}
\usepackage{amsmath}

\def\Id{{\operatorname{Id}}}
\def\tr{{\operatorname{tr}}}

\newtheorem{theorem}{Theorem}

\newtheorem{remark}[theorem]{Remark}
\newtheorem{example}[theorem]{Example}

\newtheorem{fact}[theorem]{Fact}

\begin{document}
\title{Compact Osserman manifolds with neutral metric}
\author{M. Brozos-V\'{a}zquez \and E. Garc\'{\i}a-R\'{\i}o \and P. Gilkey \and R. V\'{a}zquez-Lorenzo}
\address{MBV: E. U. Polit\'{e}cnica, Department of Mathematics, University of A Coru\~na,
Spain}
\email{mbrozos@udc.es}
\address{EGR, RVL:
Faculty of Mathematics, University of Santiago de Compostela,
Spain}
\email{eduardo.garcia.rio@usc.es, $\,$ ravazlor@edu.xunta.es}
\address{PG: Mathematics Department, University of Oregon, Eugene OR 97403
USA}
\email{gilkey@uoregon.edu}
\thanks{Supported by projects MTM2009-07756 and INCITE09 207 151 PR (Spain)}
\date{revised version of \today}
\begin{abstract}
It is shown that if a compact four-dimensional manifold with metric of neutral signature
is Jordan-Osserman, then it is either of constant sectional curvature or Ricci flat.
\end{abstract}
\maketitle

\section{Introduction}\label{s:introduction}

Self-dual Einstein metrics are of special interest both in Physics
and Geometry (see for example \cite{APV, CP1, CP2, DW} and the references therein).
From the point of view of the curvature, self-dual Einstein metrics are characterized by
the pointwise Osserman property (the eigenvalues of the Jacobi operators are constant on the unit
pseudo-spheres at each point  $p\in M$).
Four-dimensional Osserman metrics of
Riemannian signature are locally real or complex space forms, while
Lorentzian Osserman metrics are necessarily of constant curvature
(see \cite{GKV, Gilkey1} for more
information). {On the other hand}, there are plenty of
examples of Osserman metrics in signature $(--++)$ (cf. \cite{BBR,
CL-GR-G-VL, De, DR-GR-Vl}).
The Jordan normal form plays a crucial role in the higher signature
setting -- a self-adjoint linear transformation need not be
determined by its eigenvalues if the metric in question is
indefinite. One says that $M$ is \emph{Jordan-Osserman} if the
Jordan normal form of the Jacobi operators is constant on {the
pseudo-sphere bundles} $S^\pm (TM)$. Although the spacelike and timelike
Jordan-Osserman conditions are equivalent in signature $(--++)$ they do not
necessarily imply the null Jordan-Osserman condition (see
\cite{GR-G-VA-VL} for {details}).

The {relationship between local and global geometric properties is well developed in Riemannian geometry}. In
contrast, in areas such as Lorentzian or more generally pseudo-Riemannian geometry,
little is known about global properties of the geometry. In most
results, the sign of the curvature (sectional curvature, Ricci
curvature, etc.) plays a fundamental role; {the
celebrated Calabi-Markus Theorem \cite{CM} is one such example}. However, such
assumptions have no sense in the case of pseudo-Riemannian metrics
of neutral signature, since {one may
move from positive to negative curvature by reversing the sign of the original metric}. The simplest case to be
considered is that of four-dimensional $(--++)$-metrics, where some
generalizations of the Hitchin-Thorpe inequality for Einstein
metrics have been developed by Law  and Matsushita \cite{La,
La-Mat}.

{The following is one of the two main results of this paper. It shows that compactness is a strong restriction
when dealing with Osserman metrics.}

\begin{theorem}\label{main result}
Let $(M,g)$ be a compact simply connected Jordan-Osserman manifold with metric of signature
$(--++)$. Then it is of constant sectional curvature or otherwise the
Jacobi operators are two-step nilpotent.
\end{theorem}

Further observe that the Jordan normal form of the Jacobi operators may change from point
to point in an Osserman manifold. Hence, we have

\begin{theorem}\label{main result2}
Let $(M,g)$ be a compact simply connected Osserman manifold with metric of signature
$(--++)$. Then $\tau[M]=0$ and the Jacobi operators have only one
eigenvalue, which may be a single, double or triple root of the
minimal polynomial. Moreover $\chi[M]\leq 0$, and $\chi[M]=0$ if and
only if the Jacobi operators are nilpotent.
\end{theorem}

Through this paper we will construct some additional structures (almost paraHermitian and almost
anti-Hermitian ones) on Osserman manifolds by means of the conformal Weyl curvature operator.
{We shall establish Theorems \ref{main
result} and \ref{main result2} by
analyzing topological obstructions to the existence of such structures}.

\section{Osserman metrics and self-dual structures}\label{s:osserman-sd}

For any non-null vector $X$, the Jacobi operator $\mathcal{J}(X)=R(X,\,\cdot\,)X$ can be viewed as a
self-adjoint
operator on the Lorentzian vector space $\langle X\rangle^\perp$. {The} following
possibilities may occur at the algebraic level \cite{BBR}:
\begin{itemize}
\item[]\emph{Type Ia}: the Jacobi operators are diagonalizable in an orthonormal basis,
\item[]\emph{Type Ib}: the Jacobi operators have a complex eigenvalue,
\item[]\emph{Type II}: the Jacobi operators are a double root of their minimal polynomial,
\item[]\emph{Type III}: the Jacobi operators are a triple
root of their minimal polynomial.
\end{itemize}

Considering the curvature tensor $R$ as an endomorphism of
$\Lambda^2(M)$, we have the $O(2,2)$-decomposition
$
R\equiv \frac{Sc}{12}\;\Id_{\Lambda^2} +Ric_0 + W:\,\,\Lambda^2
\rightarrow \Lambda^2,
$
where $Ric_0$ denotes the traceless Ricci tensor, $Ric_0(X,Y)$
$=$ $Ric(X,Y)$ $-$ $\frac{Sc}{4}\,g(X,Y)$. The Hodge star
operator $\star:\Lambda^2\rightarrow \Lambda^2$ associated to any
$(--++)$-metric induces a further splitting $\Lambda^2
=\Lambda^2_+\oplus\Lambda^2_-$, where $\Lambda^2_\pm$ denotes the
$\pm 1$-eigenspaces of the Hodge star operator, that is,
$\Lambda^2_\pm=\{\alpha\in\Lambda^2(M)/ \star\alpha=\pm\alpha\}$.
Correspondingly, the curvature tensor further decomposes as
\begin{equation}\label{eq:spliting-2}
R\equiv \frac{Sc}{12}\;\Id_{\Lambda^2} + Ric_0 + W^++W^-,
\end{equation}
where $W^\pm=\frac{1}{2}\;(W\pm\star W)$. Recall that a
pseudo-Riemannian $4$-manifold is called \emph{self-dual} (resp.,
\emph{anti-self-dual}) if $W^-=0$ (resp., $W^+=0$).

The connection between Osserman and (anti-) self-dual metrics relies on the following Fact,
and thus the analysis of the (anti-) self-duality conditions will play a basic role in
what follows.

\begin{fact}{\rm\cite{ABBR, GKV}}\label{fact-A}
A four-dimensional pseudo-Riemannian manifold
is pointwise Osserman if and only if it is Einstein self-dual (or
anti-self-dual).
\end{fact}

Moreover, there is a one to one correspondence among
the Jordan normal forms of the Jacobi operators and the corresponding ones for $W^\pm$ as
follows.
Take a local orthonormal basis of vector fields $\{e_1,e_2,e_3,e_4\}$  on $M$. Then, local bases of the spaces of
self-dual and anti-self-dual two-forms are as follows
\[
E_1^\pm = \frac{e^1\wedge e^2 \pm e^3\wedge e^4}{\sqrt{2}},\quad
E_2^\pm = \frac{e^1\wedge e^3 \pm e^2\wedge e^4}{\sqrt{2}},\quad
E_3^\pm = \frac{e^1\wedge e^4 \mp e^2\wedge e^3}{\sqrt{2}}.
\]
Further note that $\langle
E_1^\pm, E_1^\pm \rangle =1$, $\langle E_2^\pm, E_2^\pm \rangle
=-1$, $\langle E_3^\pm, E_3^\pm \rangle =-1$, and therefore
the self-dual and anti-self-dual Weyl
curvature operators $W^\pm:\Lambda^2_\pm \longrightarrow
\Lambda^2_\pm$ may have the same Jordan canonical forms, as previously considered
for the Jacobi operators. Furthermore, assuming $W^-=0$ (a completely analogous analysis
will give the anti-self-dual case) the following algebraic correspondence, with respect to
the basis above, is obtained after a straightforward calculation.

\noindent{\underline{Type Ia}}: Diagonalizable Jacobi operators with eigenvalues
$\alpha$, $\beta$, $\gamma$ correspond
to diagonalizable self-dual Weyl curvature operator
\[
W^+=\left(\begin{array}{ccc}
2\alpha-\frac{Sc}{6}&0&0\\
0&2\beta-\frac{Sc}{6}&0\\
0&0&-2(\alpha+\beta)+\frac{Sc}{3}
\end{array}\right).
\]

\noindent{\underline{Type Ib}}: Jacobi operators have eigenvalues $\alpha$,
$\gamma\pm\beta\sqrt{-1}$
if and only if the self-dual Weyl curvature operator satisfies
\[
W^+=\left(\begin{array}{ccc}
\frac{2}{3}(\gamma-\alpha)&-2\beta&0\\
2\beta&\frac{2}{3}(\gamma-\alpha)&0\\
0&0&\frac{4}{3}(\alpha-\gamma)
\end{array}\right).
\]

\noindent{\underline{Type II}}: Jacobi operators have two eigenvalues $\alpha$ and $\beta$,
the
later being a double root of the minimal polynomial
if and only if the self-dual Weyl curvature operator satisfies
\[
W^+=\left(\begin{array}{ccc}
\frac{2}{3}(\beta-\alpha)+1 & -1 & 0\\
1 & \frac{2}{3}(\beta-\alpha)-1 & 0\\
0& 0& \frac{4}{3}(\alpha-\beta)
\end{array}\right).
\]

\noindent{\underline{Type III}}: Jacobi operators have exactly one eigenvalue $\alpha$,
which is a triple root of the minimal polynomial
if and only if the self-dual Weyl curvature operator satisfies
\[
W^+=\left(\begin{array}{ccc}
0&0&{\sqrt{2}}\\
0&0&{\sqrt{2}}\\
-{\sqrt{2}}&{\sqrt{2}}&0
\end{array}\right).
\]

\begin{remark}\label{re:Euler characteristic}\rm
Recall that a four-dimensional metric is Einstein if and only if the decomposition
(\ref{eq:spliting-2}) becomes $R\equiv \frac{Sc}{12}\;\Id_{\Lambda^2} + W^++W^-$.
In such a case, the Euler characteristic $\chi[M]$ and the Hirzebruch signature $\tau[M]$
are expressed as follows~\cite{La}
\begin{equation}\label{eq:law-formualae}
\begin{array}{l}
\displaystyle
\chi[M]=-\frac{1}{8\pi^2}\int_M\left\{\tr[(W^+)^2]+\tr[(W^-)^2]+\frac{Sc^2}{24} \right\}v,\\
\noalign{\medskip}
\displaystyle
\tau[M]=\frac{2}{3}\frac{1}{8\pi^2}\int_M\left\{ \tr[(W^+)^2]-\tr[(W^-)^2]\right \}v.
\end{array}
\end{equation}
Further observe from (\ref{eq:law-formualae}) that $\chi[M]\leq 0$ for
any compact Einstein $(--++)$--metric, provided that $W^\pm$
are not of Type Ib.
\end{remark}

\section{Proof of the Theorems}

It is a fundamental fact that an orientable 4-manifold with a
field of $2$-planes (equivalently, a two-dimensional distribution)
admits two kinds of almost complex structures { that induce opposite}
orientations. Since a neutral $4$-manifold,
i.e., an indefinite $4$-manifold of metric signature $(--++)$,
admits a field of $2$-planes, it is
necessarily an almost Hermitian $4$-manifold.

\begin{fact}{\rm\cite{Mat-Z}}\label{fact-Matsushita}
An orientable four-dimensional manifold admits a $(--++)$--metric if and
only if it satisfies a pair of Wu's conditions as follows
\[
\begin{array}{rcl}
c_{1}^{2}[M]&=&3\tau[M]+2\chi[M],\\
\noalign{\smallskip}
c_{1}^{2}[-M]&=&3\tau[-M]+2\chi[-M]=-3\tau[M]+2\chi[M],
\end{array}
\]
where $-M$ stands for $M$ with the opposite orientation.
\end{fact}

Hence, it follows from (\ref{eq:law-formualae}) the first Chern number $c_1^2[M]$ and the
first opposite Chern number
$c_1^2[-M]$ of an Einstein $(--++)$-manifold are given by
\begin{equation}\label{Chern-numbers}
\begin{array}{rcl}
c_{1}^{2}[M]&=&\displaystyle-\frac{1}{2\pi^2}\int_{M}\left\{\tr[(W^-)^2]+\frac{Sc^2}{48}\right\}v\\
\noalign{\medskip}
c_{1}^{2}[-M]&=&\displaystyle-\frac{1}{2\pi^2}\int_{M}\left\{\tr[(W^+)^2]+\frac{Sc^2}{48}\right\}v.
\end{array}
\end{equation}

\begin{remark}\label{re:estructuras}\rm
The fundamental form of an indefinite
almost Hermitian structure defines a smooth section of $\Lambda^2_-$ of constant norm $2$
and conversely, any
smooth section $\Omega$ of $\Lambda^2_-$ of constant norm $2$ is the fundamental form
of an indefinite almost Hermitian structure. {The inner product on
$\Lambda^2_\pm$ induced by a
$(--++)$-metric has signature $(--+)$. 
This shows {that there are
two almost complex structures which induce opposite
orientations};  Fact \ref{fact-Matsushita}
now follows.}
\end{remark}

\subsection{Proof of the Theorem \ref{main result}.}

Recall from Fact \ref{fact-A} that an algebraic curvature tensor is Osserman if and
only if it is Einstein and self-dual or anti-self-dual.
We assume in what follows that $W^-=0$ (a completely analogous analysis will prove the
anti-self-dual case). Hence,
the self-dual part of the Weyl conformal curvature tensor corresponds to one of the Jordan
canonical forms Ia--III.

Next, recall that $(g,P)$ is said to be an almost paraHermitian structure if $P^2=\Id$ and
$g(PX,PY)=-g(X,Y)$ for all vector fields
$X$, $Y$ on $M$. Then the fundamental form $\Omega_P(X,Y)=g(PX,Y)$ defines a section of
$\Lambda^2_+$ of constant
norm $-2$. Conversely, any smooth section $\Omega$ of $\Lambda^2_+$ of constant norm $-2$
is the fundamental form
of an almost paraHermitian structure (see also \cite{IZ}, where the inner product on
$\Lambda^2$ is taken with the opposite sign
{that we have chosen}).

Further, observe that the existence of an almost paraHermitian structure $(g,P)$ is an
equivalent condition to the existence of an almost anti-Hermitian structure $(g,J)$ (i.e.,
$J^2=-\Id$,
$g(JX,JY)=-g(X,Y)$ for all vector fields $X$, $Y$ on $M$).
Indeed, let $(g,P)$ be an almost paraHermitian structure and let $h$ be a
Riemannian metric on $M$ such that
$h(PX,PY)=h(X,Y)$ for all $X,Y$ and $h(QX,Y)=g(X,Y)$ where $Q$ is an almost product structure
(i.e., $Q^2=\Id$)
and put $J=PQ$. An
straightforward calculation shows that $J$ is an almost complex structure on $M$ and
moreover
\[
g(JX,JY)=g(PQX,PQY)=-g(QX,QY)=-g(X,Y),
\]
which shows that $(g,J)$ is almost anti-Hermitian. Moreover a straightforward calculation
shows that $J$ is
$h$-orthogonal and {that} the fundamental forms $\Omega_P$ and $\Omega_J^h(X,Y)=h(JX,Y)$ coincide
with each other. Hence,
both $P$ and $J$ induce the same orientation on $M$.

Chern classes of almost complex manifolds with anti-Hermitian metrics {were} studied in
\cite{BCGHM, BFV}, showing that
the existence of such structures is a much more restrictive condition than that of almost
Hermitian ones, as shown
in the following

\begin{fact}{\rm\cite{BCGHM}}\label{c1=0}
Let $(M,g,J)$ be an almost anti-Hermitian manifold. Then all odd Chern classes vanish
(i.e., $c_{2k+1}[M]=0$, for
all $k$).
\end{fact}

Next we consider the different possibilities for
the self-dual Weyl conformal operators.
First of all, assume $W^+$ to be of Type Ia.
A complete solution for the Osserman problem is
known in this case: either it is a space of constant sectional curvature, or {it is} an
indefinite  K\"{a}hler manifold of constant holomorphic sectional curvature or
a paraK\"{a}hler manifold of constant  paraholomorphic sectional
curvature \cite{BBR}.
Observe that in the last two cases there are exactly two-distinct
eigenvalues of the Jacobi operators in a ratio $1:\frac{1}{4}$.
If $M$ is a paracomplex space form,
then it admits an almost anti-Hermitian structure and hence $c_1[M]=0$.
Then (\ref{Chern-numbers}) shows that $Sc=0$ since $W^-=0$ and
hence $M$ is flat in the compact case.
The situation for the indefinite complex space forms  is somehow
different, since they are anti-self-dual (instead of self-dual as it occurs in the
positive  definite case) and the distinguished eigenvalue of $W^-$ has spacelike
associated eigenspace. Moreover,
\[
W^-=\frac{Sc}{12}\left(\begin{array}{ccc}2&0&0\\0&-1&0\\0&0&-1\end{array}\right)
\]
and thus
\[
\chi[M]=-\frac{1}{8\pi^2}\int_M\left\{ \tr[(W^+)^2]+\tr[(W^-)^2]+\frac{Sc^2}{24}
\right\}v=-\frac{1}{8\pi^2}\int_M \frac{Sc^2}{12}v
\]
which shows that $\chi[M]\leq 0$, with equality if and only if $M$ is flat.
K\"{a}hler-Einstein metrics have been investigated
by Petean \cite{Pe}, showing that the possible non Ricci-flat ones {occur as minimal
ruled surfaces over
curves of genus $\mathfrak{g}\geq 2$, or they occur as minimal surfaces of class
$VII_0$.}
Now, since any $VII_0$-surface has nonnegative Euler characteristic, any such surface
supports a
K\"{a}hler Osserman metric if and only if $\chi[M]=0$, and the metric is flat.
Next, assume $M$ to be a minimal ruled surface over a curve of genus $\mathfrak{g}\geq 2$.
Then the
Chern numbers $c_1^2[M]$ and $c_2[M]$ satisfy
\[
c_1^2[M]=8(1-g), \qquad
c_2[M]=4(1-g).
\]
Hence, the Hirzebruch signature $\tau[M]=\frac{1}{3}\left( c_1^2[M]-2 c_2[M] \right)$
vanishes identically,
and thus
\[
\tau[M]=\frac{1}{12\pi^2}\int_{M}\left\{
\tr[(W^+)^2]-\tr[(W^-)^2]\right\}v=-\frac{1}{12\pi^2}\int_{M}\frac{Sc^2}{24}v,
\]
which shows that $Sc=0$, and thus $M$ is flat.

Next, assume $W^+$ to be of Type Ib. Then ker$(W^++\frac{4}{3}(\alpha-\gamma) \Id)$ is
one-dimensional and timelike, since otherwise
the self-adjoint operator $W^+$ would diagonalize. The same occurs in case of Type II. In
fact, if $\alpha\neq\beta$,
then the eigenspace corresponding to $\frac{4}{3}(\alpha-\beta)$ defines an almost
paraHermitian structure on $M$.
 Next, assume $\alpha=\beta$. Then it follows from \cite{BBR} that
$\alpha=\beta=0$ and  $W^+$ is two-step nilpotent. Hence  Im$W^+$ is one-dimensional and has an
induced degenerate inner product. Therefore
the restriction of the metric to Im$W^+$ defines a one-dimensional null subspace.
Recall from Remark \ref{re:estructuras} that unit sections $\Omega^+_\pm$ of
$\Lambda^2_\pm$ of positive norm (equivalently, almost complex
structures inducing opposite orientations) exist on $M$. Hence, for any null section
$\Omega^0_\pm$ of $\Lambda^2_\pm$,
put
$\Omega^-_\pm=\frac{2}{\langle\Omega^0_\pm,\Omega^+_\pm\rangle}\Omega^0_\pm-\Omega^+_\pm$
to define a timelike section of $\Lambda^2_\pm$, and thus an almost paraHermitian
structure with fundamental
form $\Omega^-_\pm$.
Finally, if  $W^+$ is of Type III, then $(W^+)^2$ is two-step nilpotent. Hence
existence of an almost
paraHermitian structure follows then for all cases Ib, II and III.
Now, since we are assuming $W^-=0$, (\ref{Chern-numbers}) shows that
\[
c_{1}^{2}[M]=-\frac{1}{2\pi^2}\int_{M}\left\{\tr[(W^-)^2]+\frac{Sc^2}{48}\right\}v=-\frac{1}{2\pi^2}\int
_{M}\frac{Sc^2}{48}v
\]
and hence $Sc=0$, which proves that $(M,g)$ is Ricci flat. Finally, since no compact Ricci flat
Type III Jordan-Osserman manifold may exist \cite{De2}, either the sectional curvature is
constant or otherwise the Jacobi operators are two-step nilpotent.
$\hfill\square$

\subsection{Proof of Theorem \ref{main result2}.}

First of all, note that no four-dimensional Osserman metric exists whose Jacobi operators
have
three-distinct real eigenvalues, and moreover, the same holds true as concerns Osserman
metrics with
complex eigenvalues for the Jacobi operators \cite{BBR}.
{Thus we must show there are no compact Osserman metrics with two-distinct
real eigenvalues of signature $(--++)$}.
Let  $\alpha$ and $\beta$ denote the constant eigenvalues of the Jacobi operators, the
later assumed to
be of multiplicity two, and put $E_\alpha(X)=\langle X\rangle\oplus \mbox{ker}\,(\mathcal{J}(X)-\alpha\langle
X,X\rangle \Id)$,
the two dimensional subspace spanned by $X$ and the eigenvector corresponding to the
distinguished
eigenvalue $\alpha$.

Observe that, for any non-null vector $X$, the restriction of the metric tensor to
$E_\alpha(X)$ must be
{non-degenerate}, and thus either definite (of signature $(++)$ or $(--)$) or indefinite (of
Lorentzian
signature $(-+)$). Further the signature type of the $E_\alpha$'s cannot change from
definite to
indefinite, since otherwise it should pass through a degenerate case. Next, we will
consider the two
different possibilities.

If the induced metric on the $E_\alpha$'s  is definite, then ker$(\mathcal{J}(X)-\beta\langle
X,X\rangle \Id)$ is also
definite, and thus the Jacobi operators are diagonalizable. This shows that $M$ is locally
an indefinite
complex space form \cite{BBR}, and thus no such metrics exist in the compact case as shown
in the
proof of Theorem \ref{main result}. Next, assume the induced metric on the $E_\alpha$'s
is  indefinite.
Then, ker$(\mathcal{J}(X)-\alpha\langle X,X\rangle \Id)$ defines an almost paracomplex structure $P$
on
$M$, which makes
$(M,g,P)$ an almost paraHermitian manifold, and thus $c_1^2[M]=0$, which proves Ricci
flatness since $Sc=0$,
just proceeding as in the proof of Theorem \ref{main result}. Finally we show that this
cannot occur,
thus finishing the proof. Indeed, observe that the eigenvalues $\alpha$ and $\beta$ are in
a ratio $1:\frac{1}{4}$,
as it was proved in \cite{BBR} for the diagonalizable case, and moreover, it also follows
as a consequence
of the local description of the nondiagonalizable case given in \cite{DR-GR-Vl}. Hence,
$\alpha=4\beta=0$, which
is a contradiction.

As a consequence of the above, a compact and simply connected Osserman $4$-manifold $M$
with metric
of signature $(--++)$ has $\tau[M]=0$ and $\chi[M]\leq 0$, with equality if and only if
the  Jacobi operators
are nilpotent.
$\hfill\square$

\section{Compact K\"{a}hler-Osserman manifolds}

If a compact complex surface admits an indefinite K\"{a}hler-Einstein metric, then it is one of
the following \cite{Pe}:
i) a complex torus; ii) a Hyperelliptic surface; iii) a primary Kodaira surface; iv) a minimal ruled
surface over a curve of genus $\mathfrak{g}\geq 2$, or v) a minimal surface of class $VII_0$. Moreover, the
existence of the later metrics (case v)) is still an open problem.

However, a complete description of compact surfaces admitting an indefinite K\"{a}hler-Osserman metric
follows now from \cite{Kamada, Pe} and the previous results.
Indeed, since $\tau[M]=b^+-b^-=0$, it follows that
$\chi[M]=2(1-b_1+b^-)$, and thus  (since $\chi[M]\leq 0$)
$b_1\geq 1+b^-$ and $b_1= 1+b^-$ if and only if $\chi[M]=0$.
Now, an indefinite K\"{a}hler
surface satisfies
$H^2(M;\mathbb{R})\neq 0$, from where $b^-=b^+\geq 1$, and hence $b_1\geq 2$.
This shows that no minimal surface of class $VII$ is K\"{a}hler Osserman. Hence, the list of possible
K\"{a}hler-Osserman metrics reduces to i) and iii) (see the proof of Theorem \ref{main result}).

Recall that a \emph{hypersymplectic} structure on a $4n$-dimensional manifold $(M,g)$
of neutral signature is given by a triple $(J_1, J_2, J_3)$
of skew-adjoint endomorphisms of the tangent bundle satisfying
\[
J_1^2=-\Id, \qquad J_2^2=J_3^2 =\Id, \qquad J_1J_2=-J_2J_1=J_3
\]
such that the corresponding $2$-forms
$\Omega_i(X,Y)=g(J_iX,Y)$ are closed. Equivalently, the tensor fields $J_i$
are all parallel, and hence any hypersymplectic manifold is a
Ricci-flat indefinite K\"{a}hler manifold. Therefore, any hypersymplectic structure
is self-dual and Ricci flat, and thus Osserman. It was shown by Kamada
\cite{Kamada} the existence of hypersymplectic structures on complex
torus and primary Kodaira surfaces, and hence Fact \ref{fact-A} shows

\begin{theorem}
A compact complex surface admits an indefinite K\"{a}hler-Osserman metric
if and only if it is a complex torus or a primary Kodaira surface.
\end{theorem}

\begin{remark}\rm\label{re:hoy}
The special significance of hypersymplectic structures in dimension four
is given by the fact that any such structure is just a basis
$\{\Omega_1,\Omega_2,\Omega_3\}$ of closed anti-self-dual $2$-forms satisfying
$-\Omega_1^2=\Omega_2^2=\Omega_3^2$ and
$\Omega_i\wedge\Omega_j=0$ for all $i\neq j$.
Now, the endomorphisms $J_{12}=J_1+J_2$ and $J_{13}=J_1+J_3$  are skew-adjoint,
nilpotent and moreover $J_{12}$ and $J_{13}$ are parallel. Therefore Im$J_{12}$
and Im$J_{13}$ are parallel degenerate distributions on $M$, which shows that
$(M,g)$ is a Walker manifold.
Four-dimensional self-dual Ricci flat Walker metrics are known to be Riemannian
extensions of torsion-free connections with skew-symmetric Ricci tensor
\cite{walker-metrics}. Now,it follows from the analysis in \cite{De3} (see also \cite{CL-GR-G-VL})
that any hypersymplectic four-manifold is locally isometric
to the cotangent bundle $T^*\Sigma$ of an affine surface $\Sigma$ with torsion-free connection
$D$ given by
$\Gamma_{11}{}^1=-\partial_1\varphi$ and $\Gamma_{22}{}^2=\partial_2\varphi$,
for an arbitrary function $\varphi$.
The metric $g$, in induced coordinates $({ x^i},x_{i'})$ now reads as
\[
g_\nabla= 2\, dx^i\circ { dx_{i'}}+(\phi_{ij}-2x_{k'}\Gamma_{ij}{}^k)dx^i\circ dx^j
\]
for some arbitrary symmetric $(0,2)$-tensor field $\phi$ on $\Sigma$.

Further, it has been shown by Kamada \cite{Kamada} that any compact hypersymplectic surface
admits two parallel null and orthogonal vector fields (i.e., the Walker metric is strict).
A straightforward calculation now shows that any strictly Walker four-manifold is Einstein self-dual
(hence a Riemannian extension as above)
and moreover, the anti-self-dual Weyl curvature operator is two-step nilpotent.
Now, a Riemannian extension has two-step nilpotent self-dual Weyl curvature operator if and only
if the connection $D$ is equiaffine \cite{walker-metrics}, thus showing that
\emph{compact hypersymplectic manifolds are locally Riemannian extensions over flat affine surfaces}.
\end{remark}

As an application of the above we have the following examples of compact Osserman manifolds generalizing those in \cite{Pe}

\begin{example}\label{rmk-13}\rm
There are metrics on the $4$-torus which
are Walker signature $(2,2)$ metrics which are Osserman and which are
constructed as follows. Manifolds of this type were studied in this context
first in \cite{GIZ}. Let
$\{\theta_1,...,\theta_4\}$ be the usual periodic parameters on the
$4$-torus. Let $g_{ij}=g_{ij}(\theta_3,\theta_4)$ be periodic functions for
$3\le i,j\le 4$. We consider the signature $(2,2)$ Walker manifold
$\mathcal{T}:=(\mathbb{T}^4,g)$ where $g$ is given relative to the canonical
frame $\partial_{\theta_i}$ by
\begin{equation}\label{g}
g=2d\theta_1\circ d\theta_3+2d\theta_2\circ d\theta_4
+g_{33}d\theta_3\circ d\theta_3
+2g_{34}d\theta_3\circ d\theta_4
+g_{44}d\theta_4\circ d\theta_4.
\end{equation}
As the only non-zero curvature is
$R_{3443}=(2g_{34/34}-g_{33/44}-g_{44/33})/2$,
$$\mathcal{J}(\xi):\operatorname{Span}\{\partial_3,\partial_4\}
\rightarrow\operatorname{Span}\{\partial_1,\partial_2\}\rightarrow0$$
and $\mathcal{T}$ is nilpotent Osserman. We assume $R_{3443}$ does not vanish
identically so $\mathcal{T}$ is not flat. The universal cover is a
generalized plane wave manifold \cite{GN05} and hence geodesically complete;
thus
$\mathcal{T}$ is itself geodesically complete. Clearly $\mathcal{M}$ is globally
spacelike Jordan Osserman (or globally timelike Jordan Osserman) if and only if
$R_{3443}$ never vanishes.
Since $R_{3443}$ is in divergence form,
$$\int\int_{\mathbb{T}^2}R_{3443}(\theta_3,\theta_4)d\theta_3d\theta_4=0\,.$$
Thus $R_{3443}$ must change sign and $\mathcal{T}$ is not globally Jordan
Osserman. Finally observe that the metric (\ref{g}) is locally isometric to a
hypersymplectic metric as discussed in Remark \ref{re:hoy}.
\end{example}

\end{document}